\documentclass[reqno]{amsart}
\usepackage{epsfig,hyperref}

\newtheorem{conjecture}{Conjecture}

\def\sstyle{\scriptstyle}
\def\dps{\displaystyle}
\def\binom#1#2{\left(\begin{array}{@{}c@{}} #1\\#2\end{array}\right)}
\def\chordA{\begin{picture}(30,30)
\put(0,-12){\epsfxsize=30pt\epsfbox{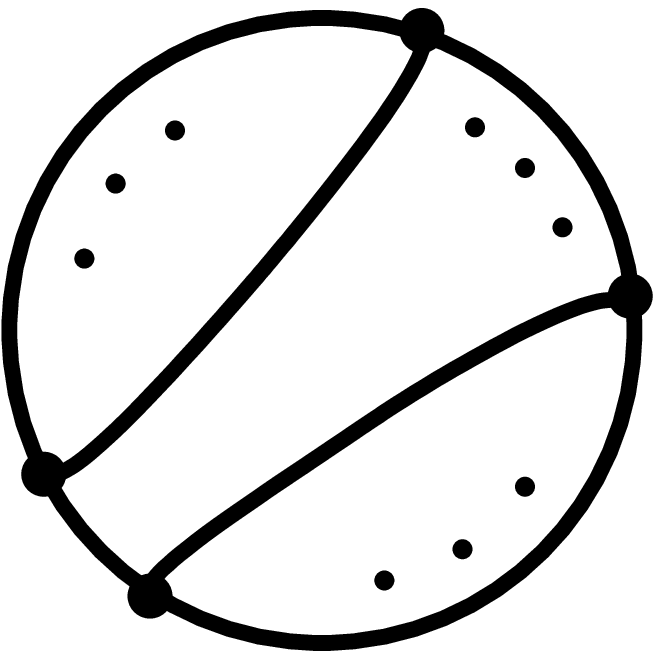}}
\put(-5,-8){$\sstyle i$}
\put(0,-18){$\sstyle i+1$}
\end{picture}}
\def\chordB{\begin{picture}(30,30)
\put(0,-12){\epsfxsize=30pt\epsfbox{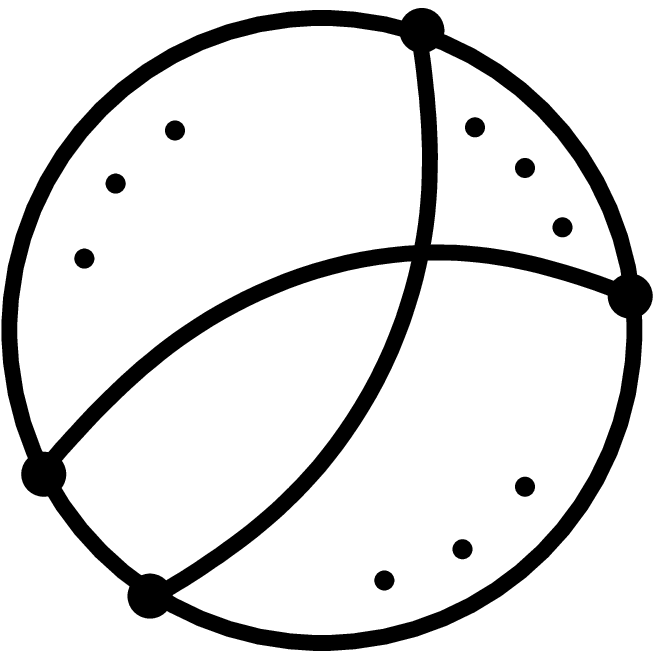}}
\put(-5,-8){$\sstyle i$}
\put(0,-18){$\sstyle i+1$}
\end{picture}}
\def\chordE{\begin{picture}(30,30)
\put(0,-12){\epsfxsize=30pt\epsfbox{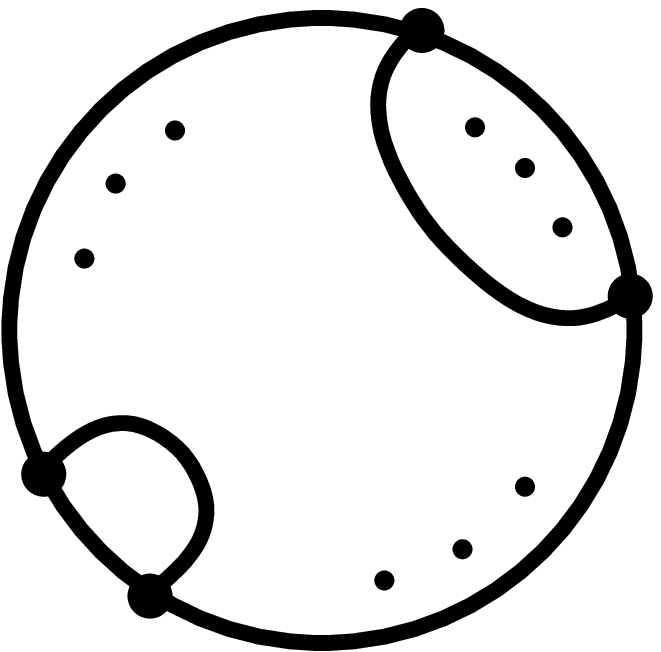}}
\put(-5,-8){$\sstyle i$}
\put(0,-18){$\sstyle i+1$}
\end{picture}}
\def\Lfoura{\begin{picture}(10,10)
\put(0,-2){\epsfxsize=10pt\epsfbox{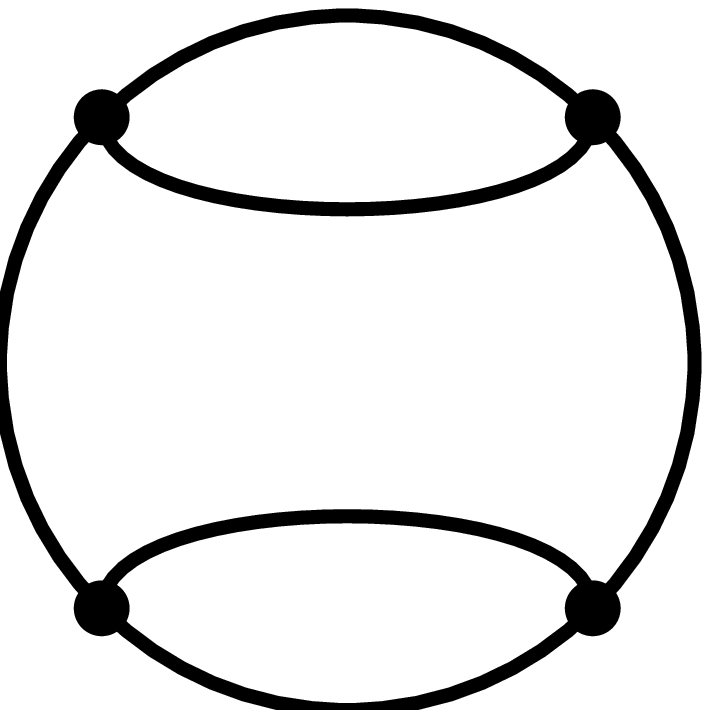}}
\end{picture}}
\def\Lfouraa{\begin{picture}(10,10)
\put(0,-2){\epsfxsize=10pt\epsfbox{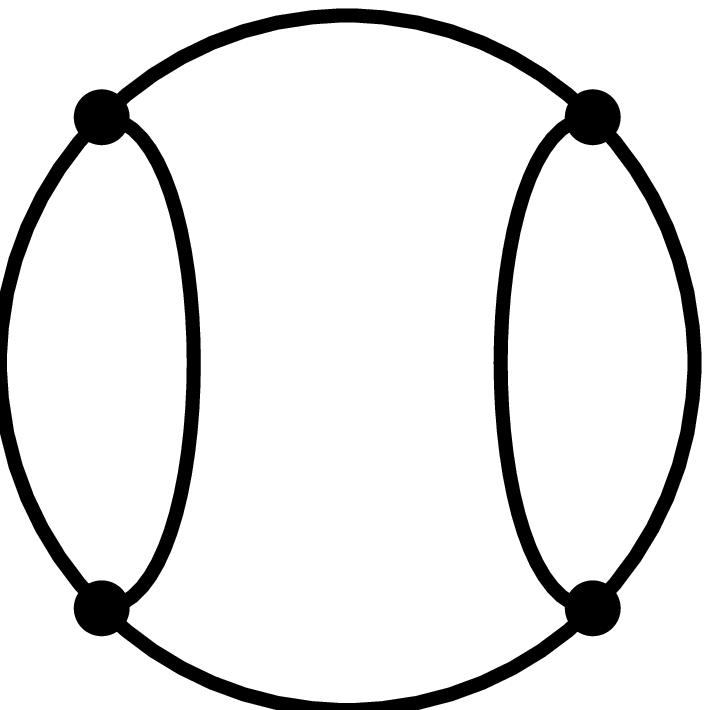}}
\end{picture}}
\def\Lfourb{\begin{picture}(10,10)
\put(0,-2){\epsfxsize=10pt\epsfbox{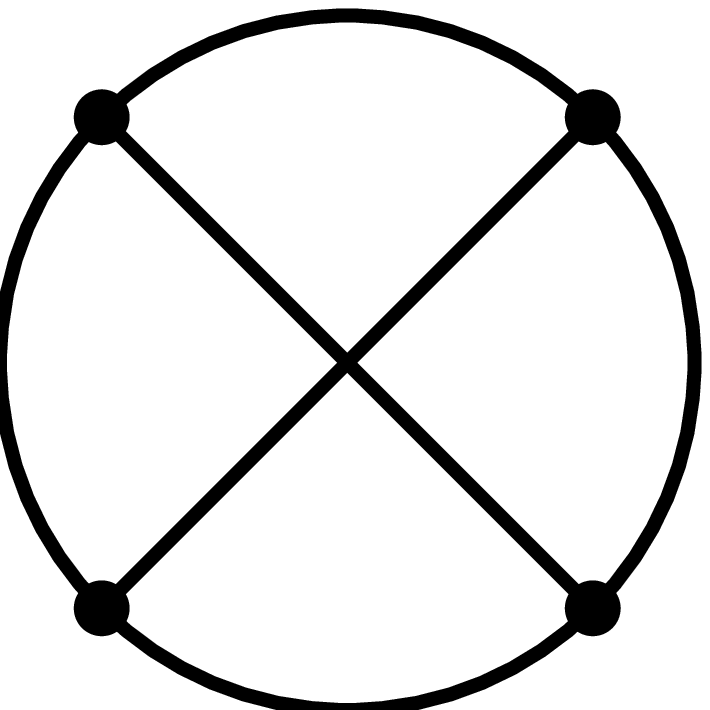}}
\end{picture}}

\begin{document}

\title{Brauer loops and the commuting variety}

\author{Jan de Gier}
\address{Department of Mathematics and Statistics, The University of
Melbourne, VIC 3010, Australia}
\email{degier@ms.unimelb.edu.au}

\author{Bernard Nienhuis}
\address{Universiteit van Amsterdam, Valckenierstraat
65-67, 1018 XE, Amsterdam, The Netherlands}
\email{nienhuis@science.uva.nl}

\date{\today}

\begin{abstract}
We observe that the degree of the commuting variety and other related
varieties occur as coefficients in the leading eigenvector of an
integrable loop model based on the Brauer algebra. 
\end{abstract}
\maketitle

\section{Introduction}

In \cite{K} Knutson introduces the upper-upper scheme $E$ which is
closely related to the commuting variety, the variety of pairs of
commuting matrices. The scheme $E$ is defined as
the variety of pairs of $n\times n$ matrices $(X,Y)$ such that
$\{(X,Y): XY$ and $YX$ upper triangular$\}$. $E$ is a reduced complete
intersection with $n!$ components $E_\pi$ labelled by permutations of
$n$. The degrees of these components provide interesting invariants
associated to permutations. The variety corresponding to the long
permutation is a degree preserving degeneration of the commuting
variety.

In unrelated research, we have encountered the degrees of $E_\pi$ as
elements of the leading eigenvector an integrable lattice model
associated to the Brauer algebra. In this note we report on this
unexpected observation, extending previous work relating the
combinatorics of alternating sign matrices to a lattice model
associated to the Temperley-Lieb algebra.

This paper is organized as follows. In the next section we introduce
the statistical mechanical Brauer loop model and formulate some
observations of its leading eigenvector in
Conjecture~\ref{con:factor}. In Section~\ref{se:scheme} we briefly
review a result by Knutson and formulate our main observation in
Conjecture~\ref{con:main}.

\section{Brauer algebra and Hamiltonian}

\subsection{Introduction}

The Temperley-Lieb algebra plays a major role in the study of solvable
models in statistical mechanics and integrable quantum chains. In
particular the XXZ quantum chain and the bond percolation model are
based on a representation of this algebra. These and some related models
have attracted renewed interest in the last few years \cite{RS1, BGN, RS2}
from an unexpected relation to the problem of counting alternating sign
matrices \cite{Bre}. This relation, which so far remains unproven, we
try to generalize in this paper.

The existing results concern the ground state of the Hamiltonian given by
\begin{equation}  
H = \sum_{i=1}^L (1-e_i) \label{TLH} 
\end{equation}
where the $e_i$ are the generators of the Temperley-Lieb algebra
satisfying
\[
e_i^2 = e_i,\qquad  e_i e_j e_i = e_i  \quad \mbox{for}\quad |i-j|=1,
\quad \mbox{and} \quad [e_i,e_j]=0 \quad \mbox{for} \quad
|i-j|>1. 
\]
While rather general boundary conditions can be considered, in this
paper we concentrate on periodic boundary conditions, requiring 
that $e_{L+i} \equiv e_i$.

\begin{figure}[h]
\centerline{
\begin{picture}(230,36)
\put(0,0){\epsfxsize=30pt\epsfbox{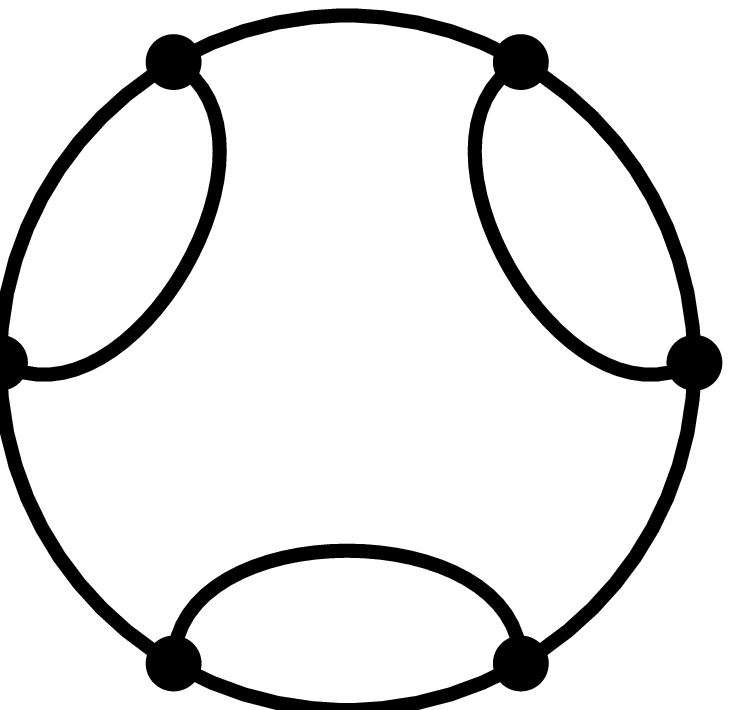}}
\put(50,0){\epsfxsize=30pt\epsfbox{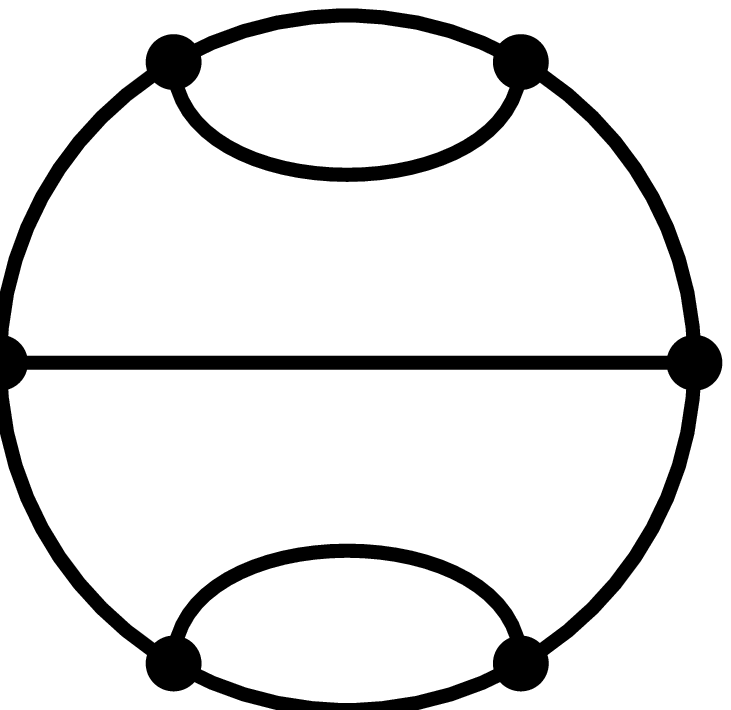}}
\put(100,0){\epsfxsize=30pt\epsfbox{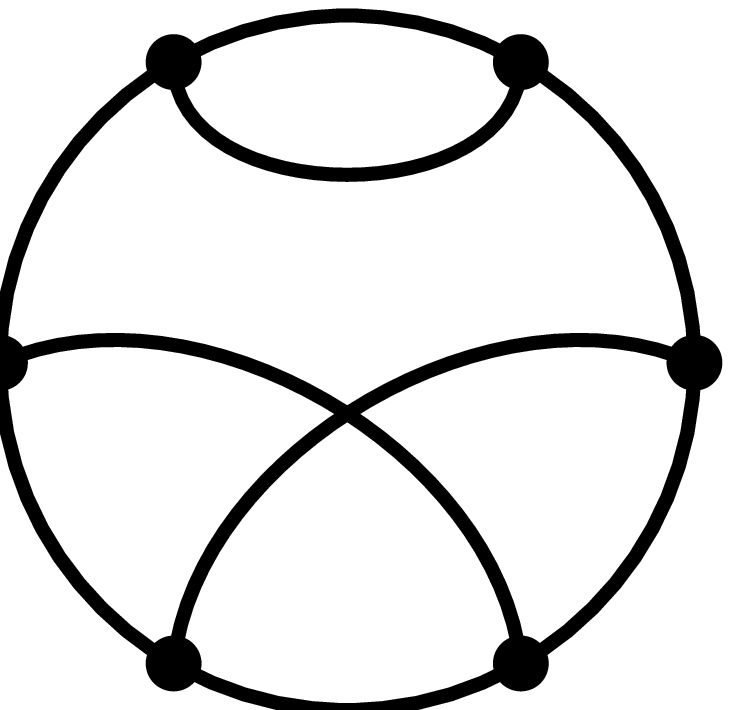}}
\put(150,0){\epsfxsize=30pt\epsfbox{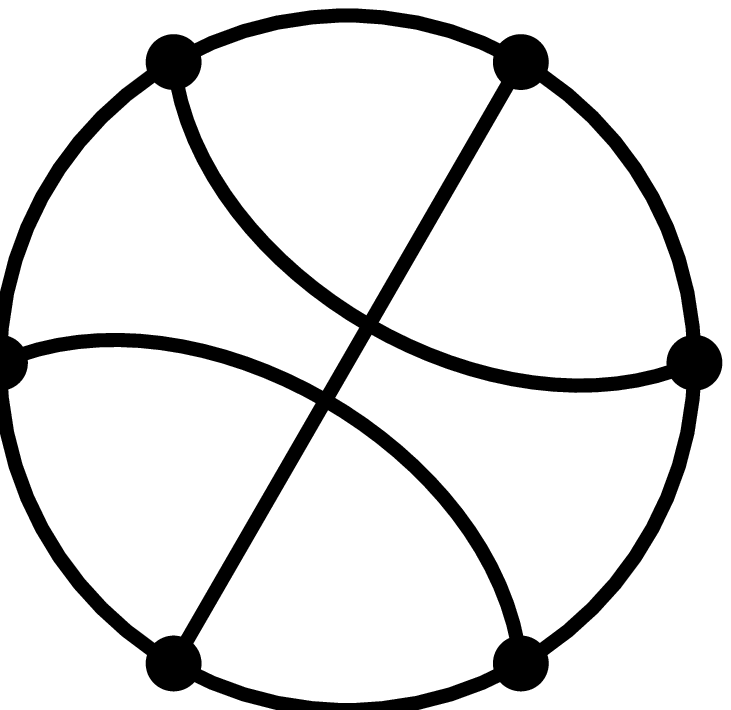}}
\put(200,0){\epsfxsize=30pt\epsfbox{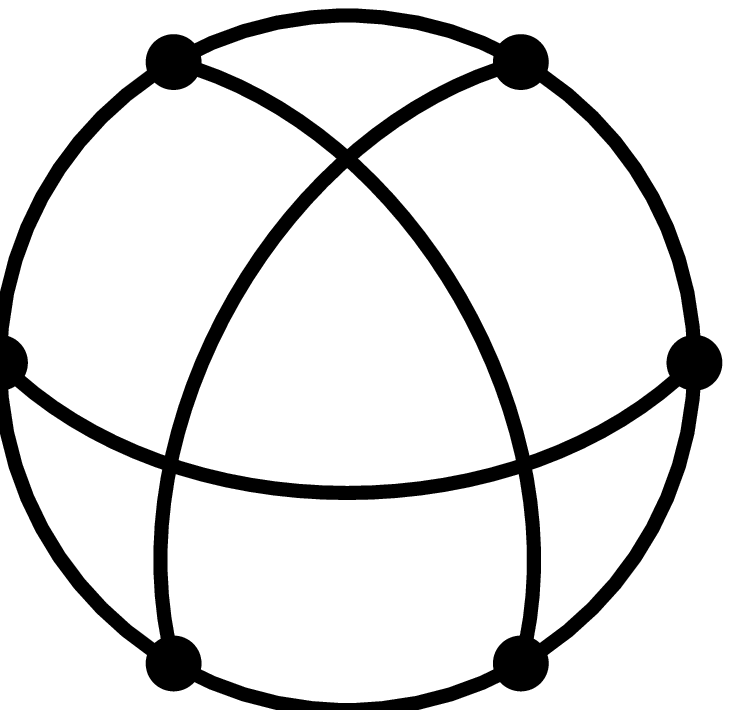}}
\end{picture}}
\caption{Chord diagrams for $L=6$.}
\label{fig:chordL6}
\end{figure}

In a convenient finite dimensional representation \cite{MS,PRGN} the $e_i$ act on (a
linear combination of) chord diagrams in which $L$ points around a
circle are pairwise connected. The five distinct chord diagrams modulo
rotations and reflections for $L=6$ are given in Figure~\ref{fig:chordL6}.

The action of $e_i$ on a chord diagram results in a new chord diagram
in which the sites $i$ and $i+1$ are connected and also the two sites
which were previously connected to $i$ and $i+1$
respectively. Pictorially,  

\vskip-14pt
\[
e_i \quad\chordA\; = \quad\chordE\; .
\]
\vskip14pt

\noindent
If $i$ and $i+1$ were already connected to each other, $e_i$ simply
acts as the identity.

The Brauer algebra \cite{Bra} is obtained by introducing besides 
the {\em monoid} $e_i$ also the {\em braid} $b_i$ which acts on the
chord diagram by interchanging the partners of $i$ and $i+1$. In terms
of a picture, 

\vskip-14pt
\[
b_i \quad\chordA\; = \quad\chordB\; .
\]
\vskip14pt

\noindent
If $i$ and $i+1$ were already paired to each other, $b_i$ leaves the chord
diagram unchanged.

The $b_i$ satisfy the usual braid relations 
\[
b_i b_j b_i = b_j b_i b_j  \; \mbox{for}\; |i-j|=1,
\quad \quad [b_i,b_j]=0 \; \mbox{for} \;
|i-j|>1,  
\]
and the additional relation
\[
b_i^2 = 1.
\]
The braids and monoids satisfy the mixed relations
\begin{eqnarray*}
&& b_i e_i = e_i b_i = e_i,\qquad  b_i b_j e_i = e_j b_i b_j = e_j e_i  
\;\; \mbox{for}\;\; |i-j|=1,\\
&& [e_i,b_j]=0 \; \mbox{for} \; |i-j|>1. 
\end{eqnarray*}

This algebra admits integrable models which have been studied in various
contexts, e.g. in an Osp(3$|$2) symmetric representation \cite{MNR} 
and to describe the low-temperature phase of O($n$) models \cite{JRS}. 
We will study the ground state of the following operator which we will
call the Brauer loop Hamiltonian
\[
H = \sum_{i=1}^{L} (3-2e_i-b_i).
\]
This operator has two properties in common with the one in (\ref{TLH}),
it is integrable and it has a smallest eigenvalue equal to zero.
We believe that it is this combination of properties that lead to
numerous special properties of the corresponding eigenvector.
 
\subsection{Groundstate}

Let us consider the action of the Brauer loop hamiltonian 
on the chord diagrams. For example,
there are three distinct chord diagrams for $L=4$, and we find
\begin{eqnarray*}
H\, \Lfoura &=& 6\; \Lfoura - 4\; \Lfouraa -2\; \Lfourb,\nonumber\\
H\, \Lfouraa &=& 6\; \Lfouraa - 4\; \Lfoura -2\; \Lfourb,\\
H\, \Lfourb &=& 12\; \Lfourb -6\; \Lfoura -6\; \Lfouraa.\nonumber
\end{eqnarray*}
From this we infer that
\[
\psi_0 = 3\; \Lfoura +3\; \Lfouraa + \Lfourb
\]
is annihilated by $H$,
\[
H \psi_0 =0.
\]
In fact, because the generators of the Brauer algebra $B_n(1)$ also
define a semi-group, the Hamiltonian $H$ is an intensity matrix
($H_{ij}\leq 0$ for $i\neq j$ 
and $\sum_i H_{ij}=0$). It therefore always has an eigenvalue $0$ with
an otherwise positive spectrum. In the following we will study the
properties of the corresponding eigenvector which we call the
groundstate.

\begin{figure}[h]
\centerline{
\begin{picture}(80,50)
\put(0,14){\epsfxsize=30pt\epsfbox{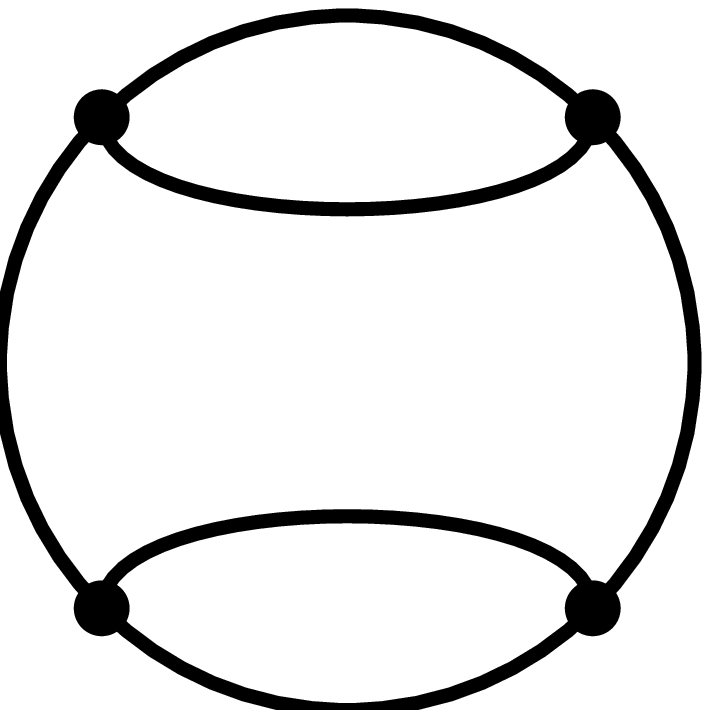}}
\put(50,14){\epsfxsize=30pt\epsfbox{L4_21.eps}}
\put(12,0){$3_2$}
\put(62,0){$1$}
\end{picture}}
\caption{Groundstate wave function for $L=4$.}
\label{fig:gsL4}
\end{figure}

We will denote states that are related by a rotation or a reflection
by the same chord diagram, and will depict the elements of the
groundstate $\psi_0$ as in Figure~\ref{fig:gsL4}, where we use
subscript to denote multiplicity. 
Let us calculate the groundstate for $L=6$ in which case there are five
distinct symmetry classes of chord diagrams (see Appendix~\ref{se:dim}). The result is given in
Figure~\ref{fig:gsL6}. We see that the smallest two
elements have the same value as those for $L=4$, and are related
to them by having an additional link that crosses each of the other links
once. In particular we note that the smallest element can be chosen to
$1$ while leaving the other elements integers (instead of rational numbers).

\begin{figure}[h]
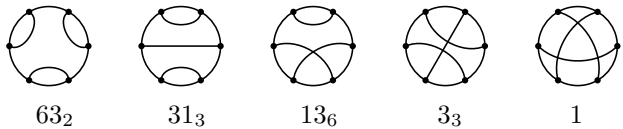

\centerline{
\begin{picture}(230,50)
\put(0,14){\epsfxsize=30pt\epsfbox{L6_214365.eps}}
\put(50,14){\epsfxsize=30pt\epsfbox{L6_321.eps}}
\put(100,14){\epsfxsize=30pt\epsfbox{L6_312.eps}}
\put(150,14){\epsfxsize=30pt\epsfbox{L6_213.eps}}
\put(200,14){\epsfxsize=30pt\epsfbox{L6_123.eps}}
\put(10,0){$63_2$}
\put(60,0){$31_3$}
\put(110,0){$13_6$}
\put(162,0){$3_3$}
\put(212,0){$1$}
\end{picture}}
\caption{Groundstate wave function for $L=6$.}
\label{fig:gsL6}
\end{figure}

\begin{figure}[hb]
\centerline{
\begin{picture}(280,190)
\put(0,14){\epsfxsize=30pt\epsfbox{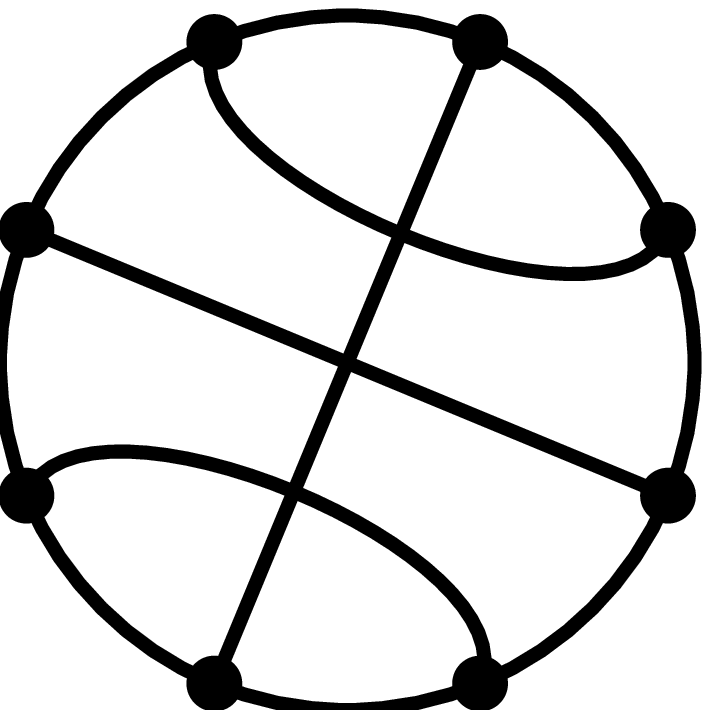}}
\put(50,14){\epsfxsize=30pt\epsfbox{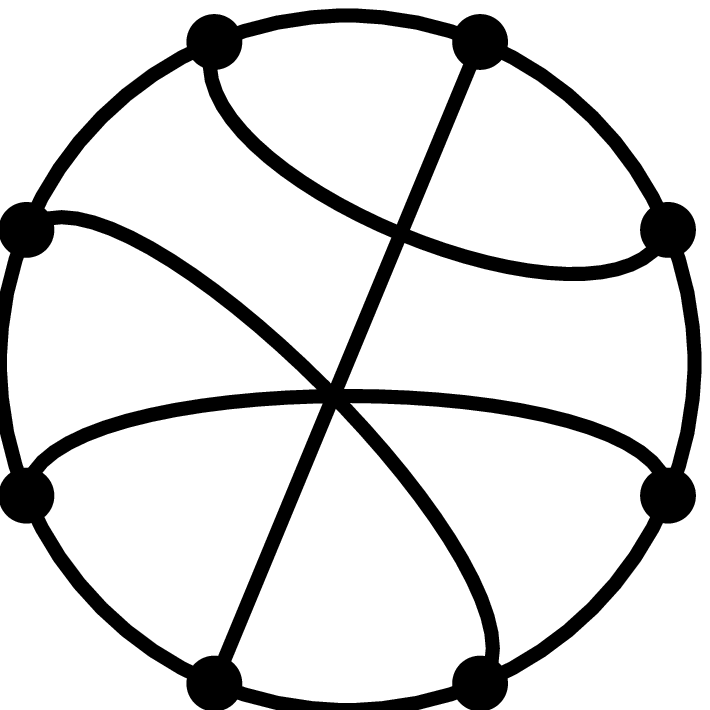}}
\put(100,14){\epsfxsize=30pt\epsfbox{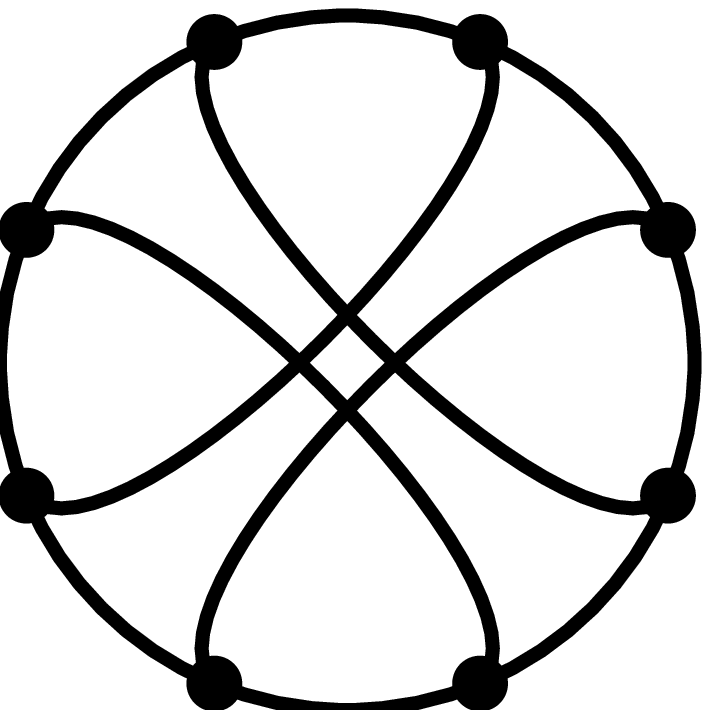}}
\put(150,14){\epsfxsize=30pt\epsfbox{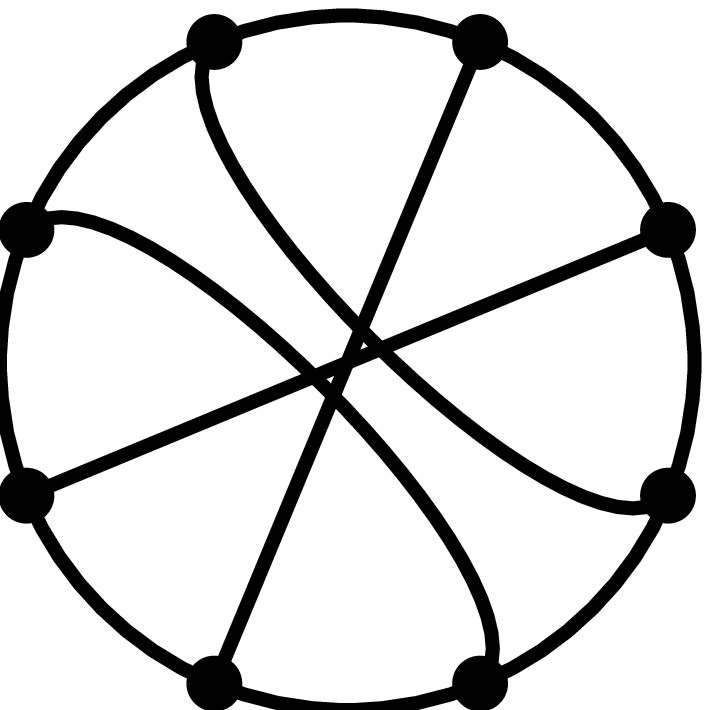}}
\put(200,14){\epsfxsize=30pt\epsfbox{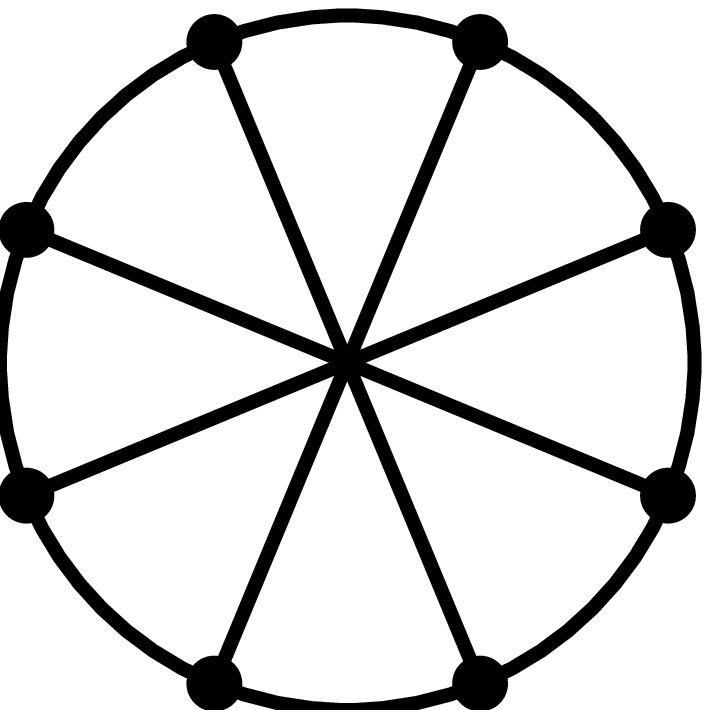}}
\put(0,84){\epsfxsize=30pt\epsfbox{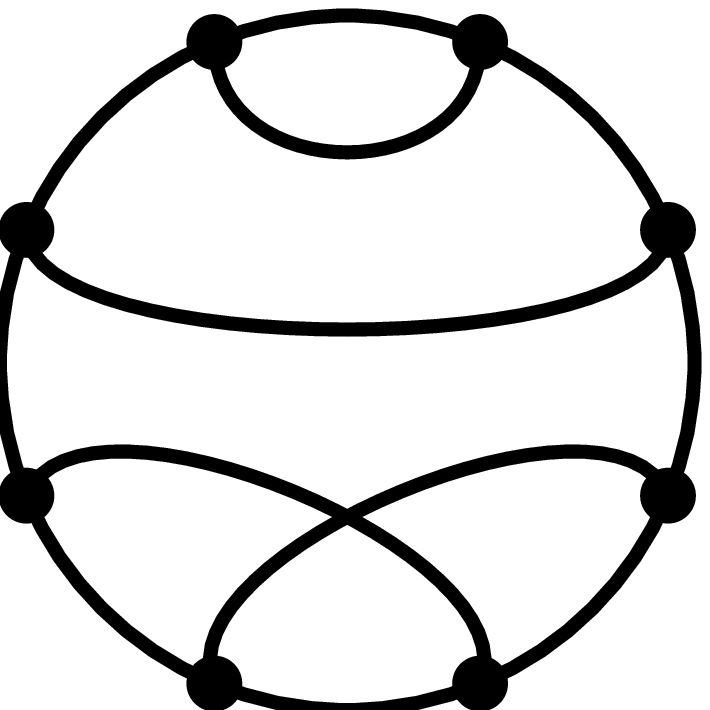}}
\put(50,84){\epsfxsize=30pt\epsfbox{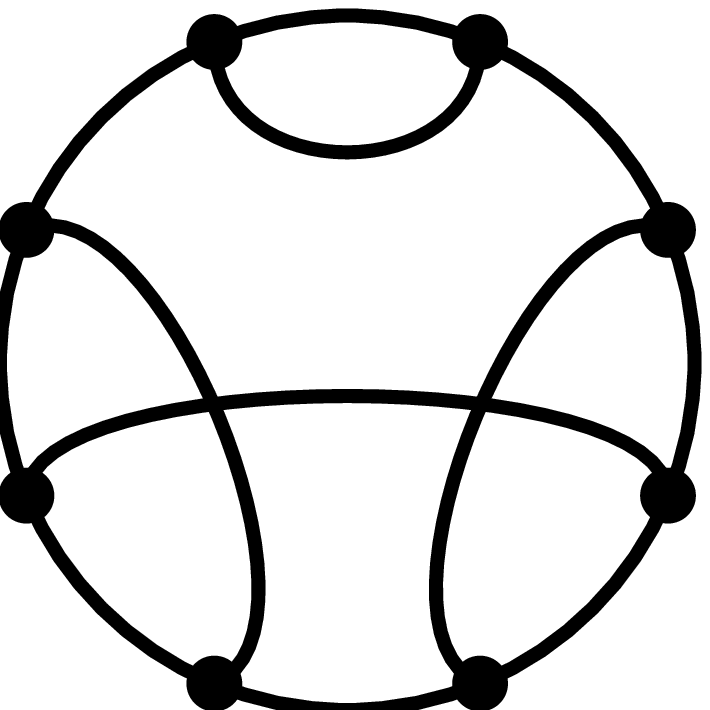}}
\put(100,84){\epsfxsize=30pt\epsfbox{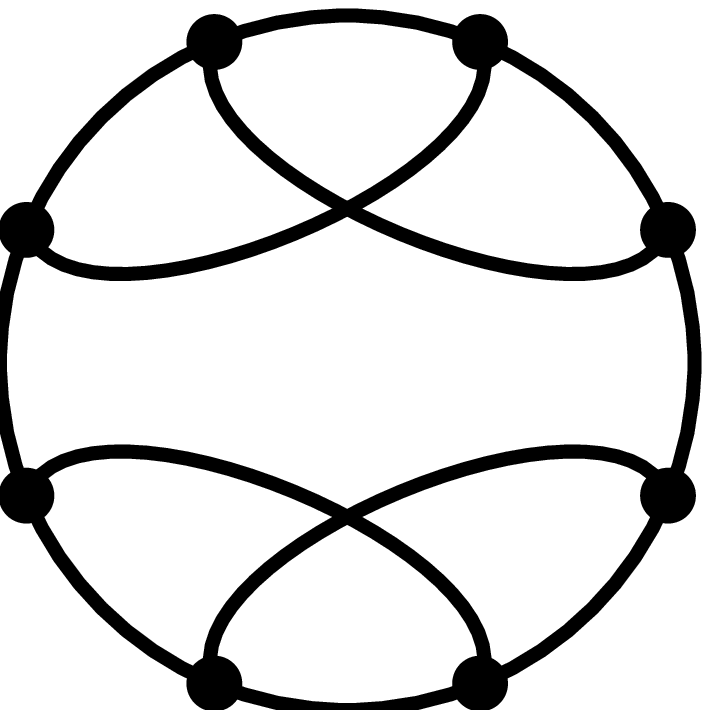}}
\put(150,84){\epsfxsize=30pt\epsfbox{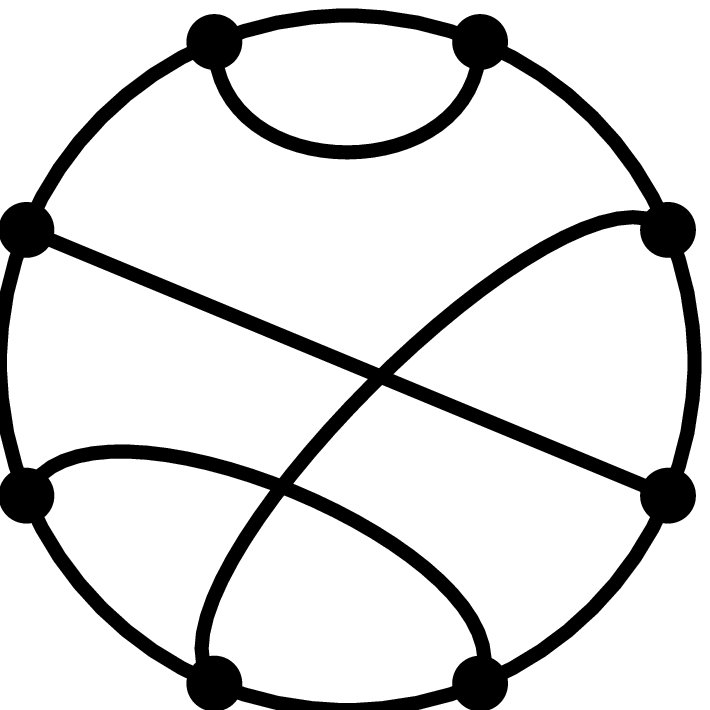}}
\put(200,84){\epsfxsize=30pt\epsfbox{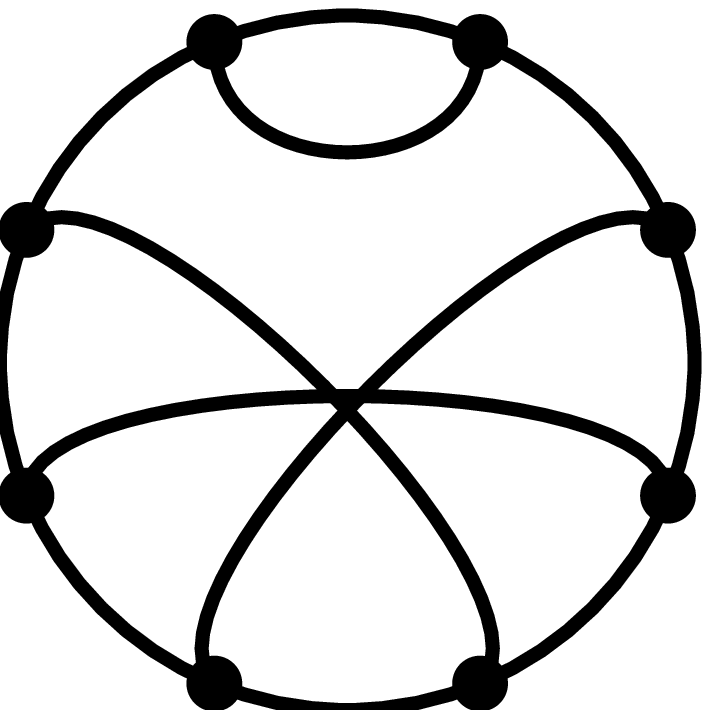}}
\put(250,84){\epsfxsize=30pt\epsfbox{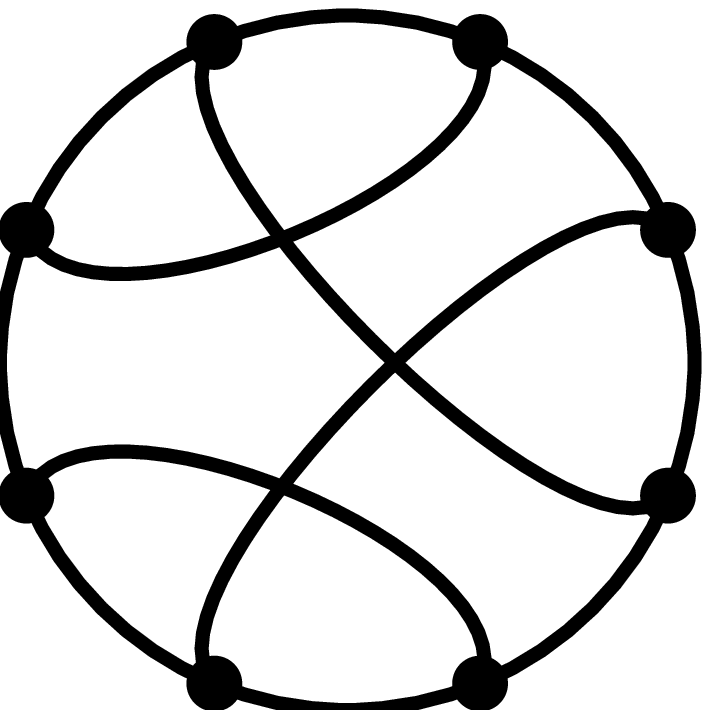}}
\put(0,154){\epsfxsize=30pt\epsfbox{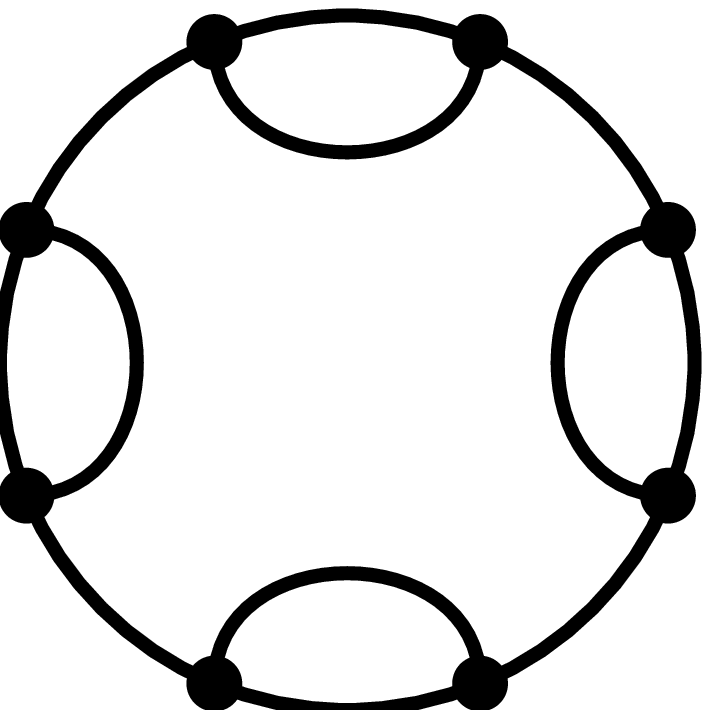}}
\put(50,154){\epsfxsize=30pt\epsfbox{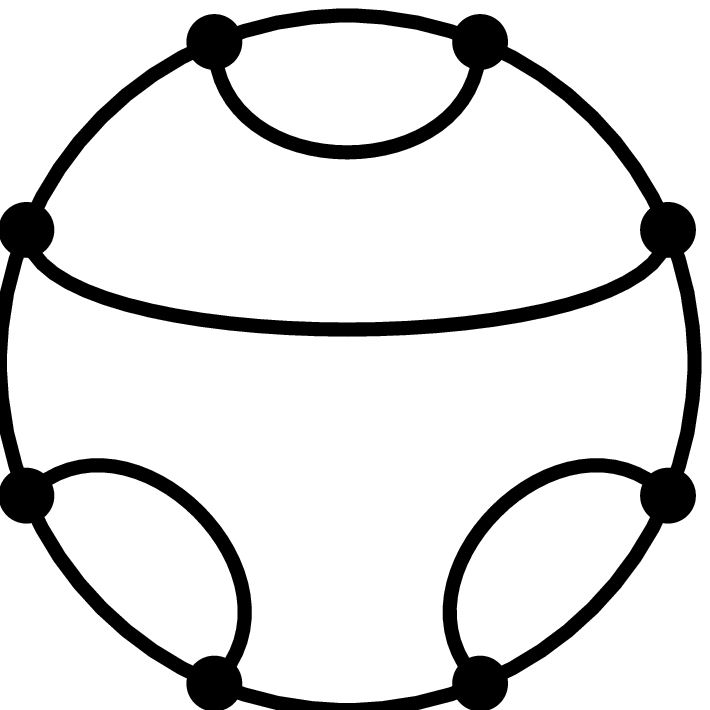}}
\put(100,154){\epsfxsize=30pt\epsfbox{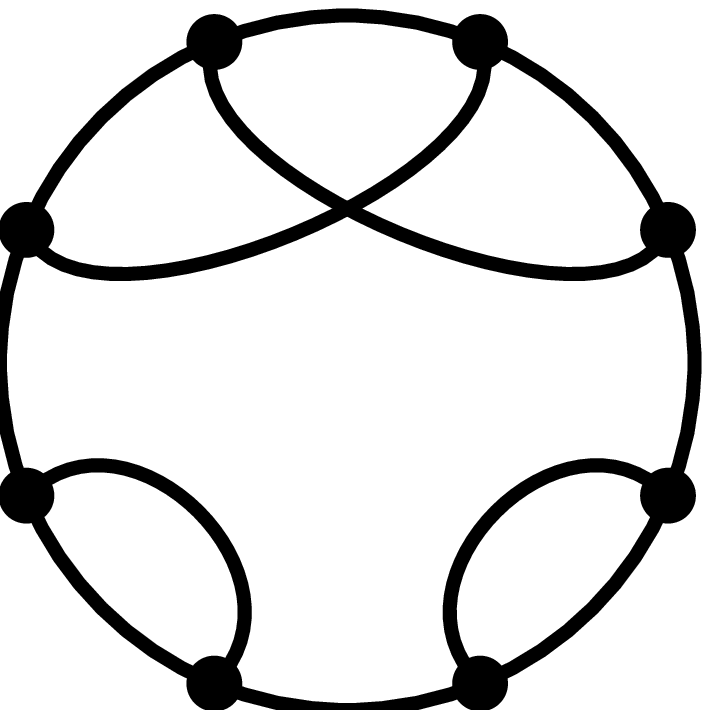}}
\put(150,154){\epsfxsize=30pt\epsfbox{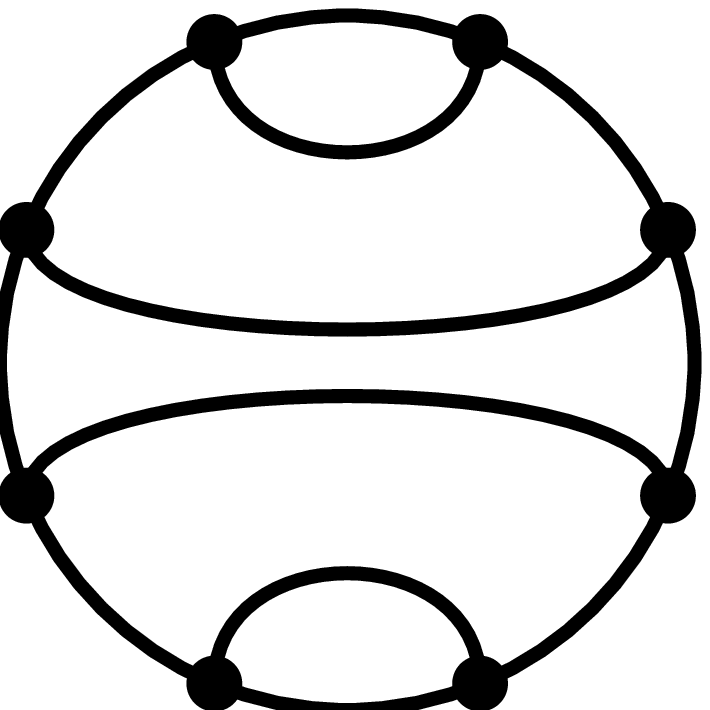}}
\put(200,154){\epsfxsize=30pt\epsfbox{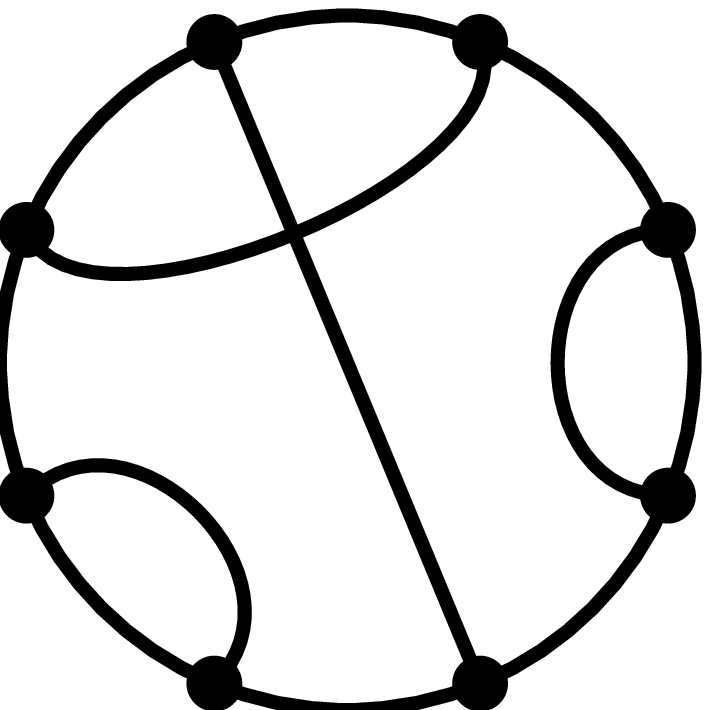}}
\put(250,154){\epsfxsize=30pt\epsfbox{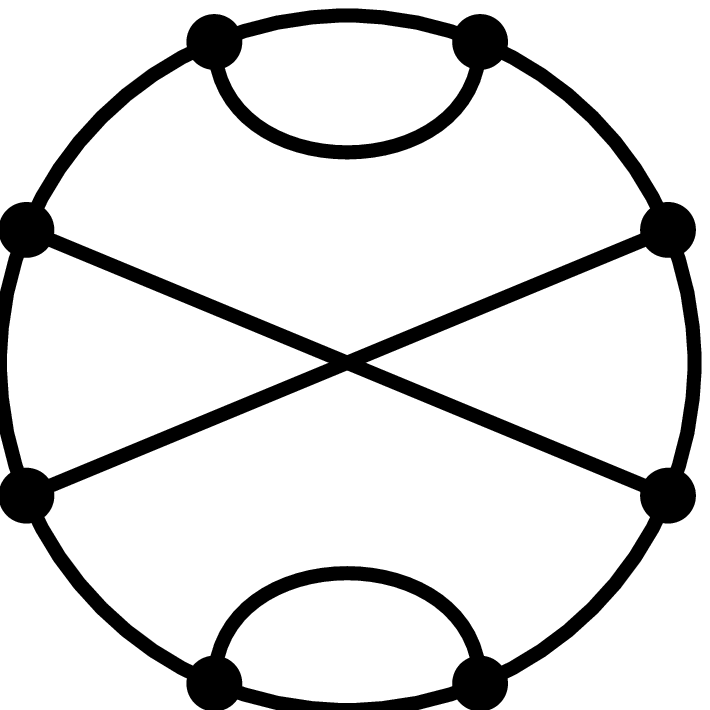}}
\put(10,0){$31_4$} \put(60,0){$13_8$} \put(112,0){$9_2$}
\put(162,0){$3_4$} \put(212,0){$1$} \put(8,70){$483_8$}
\put(58,70){$317_4$} \put(108,70){$209_4$} \put(158,70){$173_8$}
\put(210,70){$71_8$} \put(260,70){$51_8$} \put(6,140){$8297_2$}
\put(56,140){$3433_8$} \put(106,140){$1491_8$} \put(156,140){$1145_8$}
\put(206,140){$1043_4$} \put(258,140){$707_8$}
\end{picture}}
\caption{Groundstate wave function for $L=8$.}
\label{fig:gsL8}
\end{figure}

For $L=8$ there are $17$ distinct symmetry classes, and we find the
groundstate given in Figure~\ref{fig:gsL8}. As before, we note that a
chord diagram for $L=8$ has the same weight 
as a chord diagram of $L=6$ if it can be obtained by the addition of
a link that crosses each of the other links. 
Moreover the weight of a chord diagram for $L=8$ is 3 times the weight
of a chord diagram for $L=4$, if the first is obtained from the latter
by the addition of two links that cross each of the other links, but not
each other. In the next section, after we have introduced some more
notation, we will formulate a precise conjecture that includes these
observations.

\subsection{Chord diagrams related to (partial) permutations}
\begin{figure}[h]
\centerline{
\begin{picture}(60,60)
\put(10,10){\epsfxsize=40pt\epsfbox{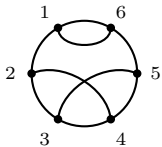}}
\put(13,51){$\sstyle 1$}
\put(0,28){$\sstyle 2$}
\put(13,3){$\sstyle 3$}
\put(42,3){$\sstyle 4$}
\put(55,28){$\sstyle 5$}
\put(42,51){$\sstyle 6$}
\end{picture}}
\caption{Chord diagram corresponding to $\pi=(312)$.}
\label{fig:312}
\end{figure}

The last four pictures in Figure~\ref{fig:gsL6} each connect the three
sites on the left hand side to those on the right hand side. We can
therefore interpret these pictures as representing elements of $S_3$,
the permutation group on three elements. For each permutation $\pi\in
S_3$, the corresponding chord diagram is the one where site $i$ is
connected to site $3+\pi(i)$. For example, $\pi=(312)$ corresponds to
the picture in Figure~\ref{fig:312}.

\begin{figure}[h]
\centerline{
\begin{picture}(130,50)
\put(0,14){\epsfxsize=30pt\epsfbox{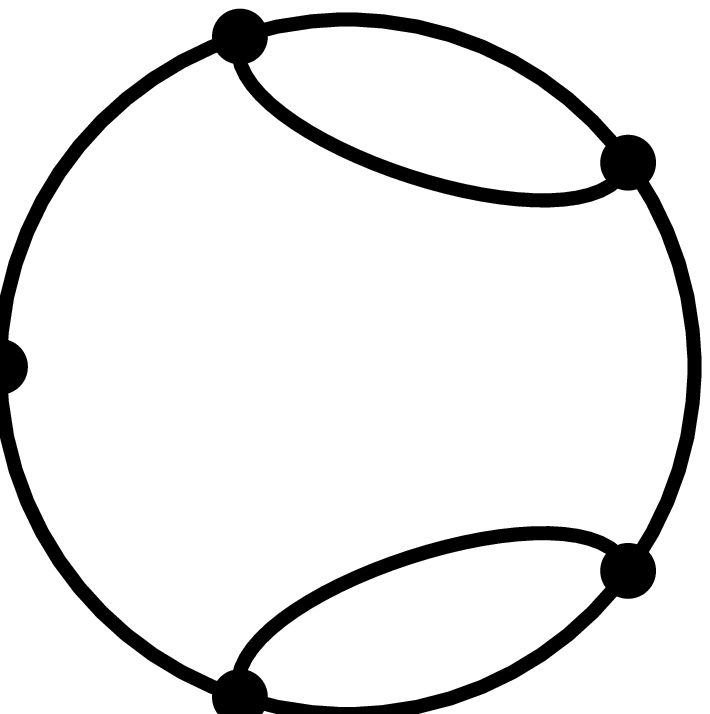}}
\put(50,14){\epsfxsize=30pt\epsfbox{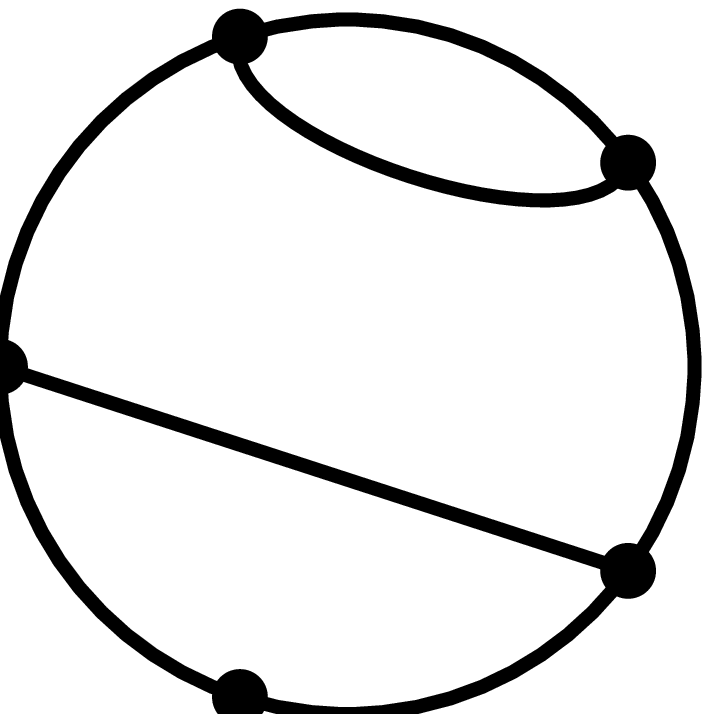}}
\put(100,14){\epsfxsize=30pt\epsfbox{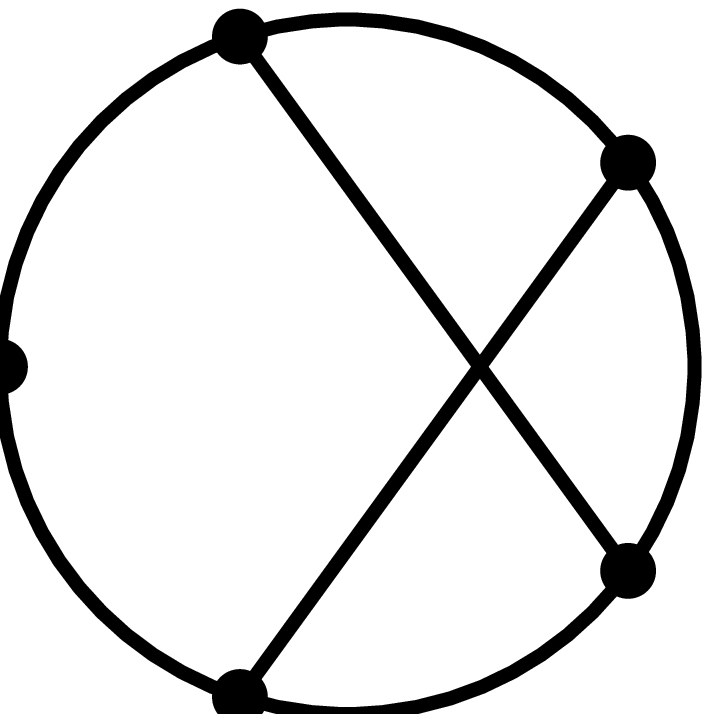}}
\put(10,0){$7_5$}
\put(60,0){$3_5$}
\put(110,0){$1_5$}
\end{picture}}
\caption{Groundstate wave function for $L=5$.}
\label{fig:gsL5}
\end{figure}

For odd systems the situation is slightly more complicated since we
cannot divide the system in two halves with equal number of
sites. However, for $L=2n+1$ we can define partial permutations of
rank $n$ for those chord diagrams where $n$ of the $n+1$ sites on the
left hand side are connected to the $n$ sites on the right hand
side. Consider for example the groundstate for $L=5$ in Figure
\ref{fig:gsL5}. 
We define for each diagram two partial permutations $\pi$ and $\pi'$. The
partial permutation $\pi$ describes the connectivity of site $i$ to
$3+\pi(i)$, and $\pi'$ describes the reverse connectivity, site $3+i$
to $\pi'(i)$ (If $\pi$ were a permutation, $\pi'$ would be the
inverse of $\pi$). For the diagrams in Figure \ref{fig:gsL5} we thus find 
\[
\begin{array}{lll}
\pi=(2\cdot1) & \pi=(21\cdot) & \pi=(1\cdot 2) \\
\pi'=(31) & \pi'=(21) & \pi'=(13)
\end{array},
\]
or in terms of matrices with a $1$ at $(i,\pi(i))$ or $(3+i,\pi'(i))$
  respectively, and zeros everywhere else (hence $P_{\pi'}=P^{\rm
  T}_\pi$),
\begin{eqnarray*}
P_\pi &=& 
\left(
\begin{array}{@{}cc@{}}
0 & 1 \\
0 & 0 \\
1 & 0
\end{array}\right),\quad
\left(
\begin{array}{@{}cc@{}}
0 & 1 \\
1 & 0 \\
0 & 0 
\end{array}\right),\quad
\left(
\begin{array}{@{}cc@{}}
1 & 0 \\
0 & 0 \\
0 & 1
\end{array}\right),\\
P_{\pi'} &=& 
\left(
\begin{array}{@{}ccc@{}}
0 & 0 & 1 \\
1 & 0 & 0 \\
\end{array}\right),\quad
\left(
\begin{array}{@{}ccc@{}}
0 & 1 & 0 \\
1 & 0 & 0 \\
\end{array}\right),\quad
\left(
\begin{array}{@{}ccc@{}}
1 & 0 & 0 \\
0 & 0 & 1 \\
\end{array}\right).
\end{eqnarray*}
Note that each chord diagram for $L=5$ corresponds to a partial
permutation. This is no longer true for larger systems. 

We will now formulate a conjecture, based on numerical evidence, about
certain elements of the groundstate $\psi_0$. 
\begin{conjecture}{\hspace{0pt}}
\label{con:factor}
\begin{itemize}
\item [i)] For each system size, the overall normalisation can be
  chosen such that the smallest element of $\psi_0$ is equal to $1$
  and all other elements are integers. 
\item [ii)] Let $\pi_n$ and $\pi_m$ be two permutations (one of them may also be a
partial permutation), and let $\pi_{n+m}$ be the permutation obtained by
their concatenation, i.e., $\pi_{n+m}(i)=\pi_n(i)$ for $1\leq i \leq n$
and $\pi_{n+m}(n+i)=n+\pi_m(i)$ for $1\leq i \leq m$. Then, with the
normalisation as in i), the weight of a chord diagram corresponding to
the (partial) permutation $\pi_{n+m}$ is equal to the product of the
weights of the chord diagrams corresponding to $\pi_n$ and $\pi_m$.
\item[iii)] Among the weights of chord diagrams that correspond to a
  permutation, the one related to the long permutation
  $w=(n,n-1,\ldots,2,1)$ is largest. 
\end{itemize}
\end{conjecture}

The weights of the chord diagrams corresponding to the long
permutation form the sequence 1, 3, 31, 1145, 154881, 77899563,
147226330175, 1053765855157617, $\ldots$ (indexed as A094579 in Sloane's
database \cite{S}). Surprisingly, the first four numbers of this
sequence form the first entries of sequence A029729, the degree of the
variety of pairs of commuting $n\times n$ matrices. The computation of
these degrees is difficult and at the time of this writing the four
terms in A029729 are the only ones known \cite{W}. Assuming the
connection, the Brauer loop model provides us with a way of obtaining
these degrees for relatively large $n$. This motivated us to find an
interpretation for other entries of the groundstate as well, and in the
next section we will formulate our main result, interpreting many more
entries of $\psi_0$ as degrees of certain algebraic varieties related
to the commuting variety.  

\section{The scheme $\overline{E}_\pi$}
\label{se:scheme}

The commuting variety is the variety of pairs of $n\times n$
matrices $(X,Y)$ such that $XY=YX$. In \cite{K} Knutson
introduces generalizations of the commuting variety: the diagonal
commutator scheme and its flat degeneration, the upper-upper scheme. The diagonal
commutator scheme is the variety of pairs of matrices $(X,Y)$ such
that $XY-YX$ is diagonal and the upper-upper scheme $E$ is
$\{(X,Y): XY$ and $YX$ upper triangular$\}$. $E$ can be
generalized to rectangular $m\times n$ and $n\times m$ matrices, and
can be naturally decomposed into subsets $E_{\pi}$, labelled by 
partial permutations $\pi$ of rank $m$ if $m<n$. The degrees of these
varieties provide interesting invariants. In \cite{K} Knutson
conjectures that the variety $\overline{E}_\pi$ for $\pi\in S_n$ is
defined as a scheme by the following three sets of equations. 

\begin{conjecture}[Knutson]
\label{con:knut}
The variety $\overline{E}_\pi$ is defined as a scheme by three
sets of equations,
\begin{itemize}
\item $XY$ and $YX$ upper triangular
 \item ${\rm diag}(XY) = {\rm diag}(P_\pi YX P^{\rm T}_{\pi})$
\item those defining the $P_\pi, P^{\rm T}_{\pi}$ matrix Schubert
varieties: for each pair $i,j$ the rank of the lower left $i \times
j$ rectangle in $X$ (resp.in $Y$) is bounded above by the number
of $1$s in that rectangle in $P_\pi$ (resp. in $P^{\rm T}_{\pi}$).
\end{itemize}
\end{conjecture}

Knutson furthermore calculates the degrees of $\overline{E}_\pi$ for
$\pi\in S_3$ \cite[Prop. 3]{K}, and finds
\begin{eqnarray*}
&&d_{(123)} = 1,\quad d_{(132)} = d_{(213)} = 3,\\
&&d_{(231)} = d_{(312)} = 13,\quad d_{(321)}=31.
\end{eqnarray*}
These degrees are exactly the same as the groundstate wavefunction
elements for $L=6$ corresponding to the permutation $\pi$, as given in
Figure~\ref{fig:gsL6}. We have calculated the degrees of $\overline{E}_\pi$ for
$\pi\in S_4$ using Macaulay2, see the appendix, and indeed found a
correspondence with the groundstate for $L=8$. We have furthermore
done similar calculations for odd system sizes, and are led to the
following intriguing conjecture.

\begin{conjecture}
\label{con:main}
The groundstate element of the Brauer Hamiltonian for $L=2n$
corresponding to a permutation $\pi\in S_n$ is equal to the degree
of the subset $\overline{E}_\pi$ of the upper-upper scheme for pairs of
$n\times n$ matrices. For $L=2n+1$ we find that the groundstate
elements corresponding to the partial permutation matrix $\pi$ of
rank $n$, corresponds to the degree of the subset $\overline{E}_\pi$
of the upper-upper scheme for an $n\times (n+1)$ and an $(n+1)\times
n$ matrix.
\end{conjecture}

For the long permutation $w=(n,n-1,\ldots,2,1)$, the last of the three
sets of defining equations in Conjecture \ref{con:knut} is empty and
the remaining two are those defining a degeneration of the commuting variety
\cite{K}, consistent with our observation in the previous
section. It was furthermore proved by Knutson, \cite[Prop. 3]{K},
that the sum over all degrees of $\overline{E}_{\pi}$ has a simple
expression,  
\[
\sum_{\pi\in S_n} \deg \overline{E}_{\pi} = 2^{n^2-n}.
\]
Assuming Conjecture \ref{con:main}, the same holds for the
corresponding sum over groundstate elements for $L=2n$, an observation
made by independent numerical calculations by Zuber \cite{Z}. For
$L=2n+1$ we find that the sum over all groundstate elements
corresponding to a partial permutation of rank $n$ is equal to
$2^{n^2}$, which is indeed what one would expect for the total degree
of $E$ for an $n\times (n+1)$ and an $(n+1)\times n$ matrix following
Knutson's argument.  

\section{Acknowledgment}

Our warm thanks go to Michel Bauer, Philippe Di Francesco and
Jean-Bernard Zuber for fruitful discussions, and to Nolan Wallach and
Allen Knutson for useful correspondence. This research was financially
supported by the Australia Reserach Council (JdG) and Stichting FOM (BN).
 
\appendix
\section{A Macaulay2 session}
In this section we describe a Macaulay2 \cite{GS} session calculating the
degree of the variety $\overline{E}_{(2,4,3,1)}$. We start with defining a
polynomial ring of $32$ variables, and define two matrices $X$ and
$Y$ and their products.\\

{\tt
\renewcommand{\arraystretch}{1.4}
\begin{tabular}{l}
i1\ :\ R = ZZ[x\_1..x\_16, y\_1..y\_16]; \\
i2\ :\ X = genericMatrix(R, x\_1, 4, 4) \\
o2\ =\hspace{-10pt}
\renewcommand{\arraystretch}{0.8}
\begin{tabular}[t]{@{\hspace{7pt}}l@{\hspace{7pt}}l@{\hspace{7pt}}l@{\hspace{7pt}}l@{\hspace{7pt}}l@{\hspace{7pt}}l@{}}
| & x\_1 & x\_5 & x\_9 & x\_13 & |\\
| & x\_2 & x\_6 & x\_10 & x\_14 & |\\
| & x\_3 & x\_7 & x\_11 & x\_15 & |\\
| & x\_4 & x\_8 & x\_12 & x\_16 & |
\end{tabular}
\renewcommand{\arraystretch}{1.4}
\\
o2\ :\ Matrix R$^4$ <--- R$^4$ \\
i3\ :\ Y = genericMatrix(R, y\_1, 4, 4) \\
o3\ =\hspace{-10pt}
\renewcommand{\arraystretch}{0.8}
\begin{tabular}[t]{@{\hspace{7pt}}l@{\hspace{7pt}}l@{\hspace{7pt}}l@{\hspace{7pt}}l@{\hspace{7pt}}l@{\hspace{7pt}}l@{}}
| & y\_1 & y\_5 & y\_9 & y\_13 & |\\
| & y\_2 & y\_6 & y\_10 & y\_14 & |\\
| & y\_3 & y\_7 & y\_11 & y\_15 & |\\
| & y\_4 & y\_8 & y\_12 & y\_16 & |
\end{tabular}
\renewcommand{\arraystretch}{1.4}
\\
o3\ :\ Matrix R$^4$ <--- R$^4$ \\
i4\ :\ XY=X*Y; \\
o4\ :\ Matrix R$^4$ <--- R$^4$ \\
i5\ :\ YX=Y*X; \\
o5\ :\ Matrix R$^4$ <--- R$^4$ \\
\end{tabular}}\\

\noindent Next we take the elements of $XY$ and $YX$ below the
diagonal and concatenate them.\\

{\tt
\begin{tabular}{l}
i6\ :\ XYupperTri = (flatten XY)\_\{1,2,3,6,7,11\};\\
o6\ :\ Matrix R$^1$ <--- R$^6$ \\
i7\ :\ YXupperTri = (flatten YX)\_\{1,2,3,6,7,11\};\\
o7\ :\ Matrix R$^1$ <--- R$^6$ \\
i8\ :\ upperTri = XYupperTri|YXupperTri;\\
o8\ :\ Matrix R$^1$ <--- R$^{12}$ \\
\end{tabular}}\\

\noindent The second set of polynomials that are part of the ideal
is obtained by taking the diagonal elements of $XY-P_\pi
YXP^{\rm T}_\pi$, where $P_\pi$ is the permutation matrix
corresponding to $\pi=(2,4,3,1)$ with a 1 at $(i,\pi(i))$ and zeros everywhere else.\\

{\tt
\begin{tabular}{l}
i9\ :\ Perm = \{2,4,3,1\};\\
i10\ :\ diagXY = map((ring XY)\^{}(numgens source XY), source XY,\\
\hspace{15mm} (i,j) -> if i
=== j then XY\_(i,i) else 0); \\
o10\ :\ Matrix R$^4$ <--- R$^4$ \\
i11\ :\ diagPYX = map((ring YX)\^{}(numgens source YX), source YX,
\\
\hspace{15mm} (i,j) -> if i === j then YX\_(Perm\#i-1,Perm\#i-1) else 0); \\
o11\ :\ Matrix R$^4$ <--- R$^4$ \\
i12\ :\ diagXYPYX = compress flatten (diagXY-diagPYX)\\
o12\ :\ Matrix R$^1$ <--- R$^4$ \\
\end{tabular}}\\

\noindent Lastly we create the set of polynomials forming the
ideals of the Schubert varieties of $P_\pi$ and $P^{\rm T}_{\pi}$. All
equations are concatenated and the ideal $I$ corresponding to $\overline{E}_\pi$
is defined.
\\

 {\tt
\begin{tabular}{l}
i13\ :\ Schubert = matrix \{\{y\_4, y\_3, det(submatrix(X,(2,3),(0,1))),\\
\hspace{15mm} det(submatrix(X,(1,3),(0,1))), det(submatrix(X,(1,2),(0,1)))\}\}; \\
o13\ :\ Matrix R$^1$ <--- R$^5$ \\
i14\ :\ total = upperTri|diagXYPYX; \\
o14\ :\ Matrix R$^1$ <--- R$^{16}$ \\
i15\ :\ total = total|Schubert;\\
o15\ :\ Matrix R$^1$ <--- R$^{21}$ \\
i16\ :\ I = ideal total;\\
o16\ :\ Ideal of R\\
i17\ :\ degree I\\
o17\ :\ 173
\end{tabular}
}

\section{Dimension of the Hilbert space}
\label{se:dim}

The total number $c_n=1,2,5,17,79,\ldots$, of symmetry classes of
chord diagrams, or pairings on a bracelet, is given by \cite{L,S}
\[
c_n = \frac14 \left(\frac1n \sum_{pq=2n} \alpha(p,q) \phi(q) + d_n+d_{n-1}\right),
\]
where
\begin{eqnarray*}
\alpha(p,q) &=&
\left\{\begin{array}{ll}
\dps \sum_{k=0}^{\lfloor p/2\rfloor} \binom{p}{2k} q^k(2k-1)!! & q\;{\rm
  even}\\[5mm]
\dps q^{p/2}(p-1)!! & q\;{\rm odd},
\end{array}
\right.\\
d_n &=& \sum_{k=0}^{\lfloor n/2\rfloor} \frac{n!}{(n-2k)!k!},
\end{eqnarray*}
and $\phi(n)$ is Euler's totient function.

\end{document}